\documentclass[graybox]{svmult}

\usepackage{mathptmx}       
\usepackage{helvet}         
\usepackage{courier}        
\usepackage{type1cm}        
\usepackage{makeidx}         
\usepackage{graphicx}        
\usepackage{multicol}        
\usepackage[bottom]{footmisc}

\makeindex             

\newif\iffinal

\finaltrue

\usepackage{todonotes}
\usepackage[acronym]{glossaries}

\newacronym{PIF}{PIF}{persistence indicator function}

\usepackage[T1]{fontenc}

\usepackage{siunitx}

\usepackage[caption=false]{subfig}

\usepackage{dsfont}

\usepackage{amsmath}
\usepackage{amssymb}

\DeclareMathOperator {\dist}       {dist}

\DeclareMathOperator {\pers}       {pers}

\DeclareMathOperator {\volume}     {vol}
\DeclareMathOperator {\wasserstein}{W\!}

\newcommand{\diagram}         {\mathcal{D}}

\newcommand{\empiricalsample} {\ensuremath{\mathds{P}}}
\newcommand{\landau}[1]       {\ensuremath{\mathcal{O}\left(#1\right)}}

\newcommand{\Lp}[1][p]        {\ensuremath{\mathrm{L}^{\!#1}}}

\newcommand{\pif}[1][\diagram]{\mathds{1}_{#1}}
\newcommand{\real}            {\ensuremath{\mathds{R}}}

\renewcommand{\d}[1]{\ensuremath{\operatorname{d}\!{#1}}}

\usepackage{tikz}
\usepackage{pgfplots}

\pgfplotsset{compat=1.14}

\usepgfplotslibrary{colorbrewer}
\usepgfplotslibrary{external}
\usepgfplotslibrary{groupplots}
\usepgfplotslibrary{fillbetween}

\tikzsetexternalprefix{Figures/External/}

\definecolor{amber}     {RGB}{255,191,  0}
\definecolor{cardinal}  {RGB}{196, 30, 58}
\definecolor{yale}      {RGB}{ 70,130,180}

\pgfplotscreateplotcyclelist{color list}{%
  cardinal,
  yale,
  amber
}

\pgfplotsset{cycle list name=color list}

\usetikzlibrary{calc}
\usetikzlibrary{shapes}

\usepackage{booktabs}
\usepackage{xspace}

\renewcommand{\th}{\textsuperscript{\textup{th}}\xspace}

\newcommand{\comment}[1]{}
\usepackage[absolute]{textpos}

\begin{document}

\title*{Topological Machine Learning with Persistence Indicator Functions}

\author{Bastian Rieck \and Filip Sadlo \and Heike Leitte}

\institute{Bastian Rieck \and Heike Leitte \at TU Kaiserslautern, \email{\{rieck, leitte\}@cs.uni-kl.de}%
  \and Filip Sadlo \at Heidelberg University, \email{sadlo@uni-heidelberg.de}%
}

\maketitle

\abstract{%
  Techniques from computational topology, in particular persistent homology, are becoming
  increasingly relevant for data analysis. Their stable metrics permit the use of many
  distance-based data analysis methods, such as multidimensional scaling, while providing a firm
  theoretical ground. Many modern machine learning algorithms, however, are based on kernels.  This
  paper presents \emph{persistence indicator functions}~(PIFs), which summarize persistence
  diagrams, i.e., feature descriptors in topological data analysis.
  PIFs can be calculated and compared in linear time and have many beneficial properties, such as
  the availability of a kernel-based similarity measure. We demonstrate their usage in common data
  analysis scenarios, such as confidence set estimation and classification of complex structured
  data.
}

\iffinal
\else
  \begin{textblock*}{\paperwidth}[1, 0](\paperwidth,0.5cm)
  \scriptsize%
  Authors' copy. Please refer to \emph{Topological Methods in Data Analysis and Visualization
  V: Theory, Algorithms, and Applications} for the definitive version of this chapter.
  \end{textblock*}
\fi

\section{Introduction}

Persistent homology~\cite{Edelsbrunner10,Edelsbrunner02,Edelsbrunner14}, now over a decade old, has
proven highly relevant in data analysis. The last years showed that the usage of topological
features can lead to an increase in, e.g., classification performance of machine learning
algorithms~\cite{Reininghaus15}.
The central element for data analysis based on persistent homology is the \emph{persistence
diagram}, a data structure that essentially stores related critical points~(such as minima or
maxima) of a function, while providing two stable metrics, namely the \emph{bottleneck distance} and
the $p$\th \emph{Wasserstein distance}.
Certain stability theorems~\cite{Cohen-Steiner07, Cohen-Steiner10} guarantee that the calculations
are robust against perturbations and the inevitable occurrence of noise in real-world data.

This stability comes at the price of a very high runtime for distance calculations between
persistence diagrams: both metrics have a complexity of at least~$\landau{n^{2.5}}$, or, if naively
implemented, $\landau{n^3}$~\cite[p.\ 196]{Edelsbrunner10}. Using randomized algorithms, it is
possible to achieve a complexity of~$\landau{n^\omega}$, where $\omega < 2.38$ denotes the best
matrix multiplication time~\cite{Mucha04}. Further reductions in runtime complexity are possible if
approximations to the correct value of the metric are permitted~\cite{Kerber16}.
Nevertheless, these algorithms are hard to implement and their performance is worse than
$\landau{n^2}$, meaning that they are not necessarily suitable for comparing larger sets of
persistence diagrams.

In this paper, we describe a summarizing function for persistence diagrams, the \gls{PIF}.
\glspl{PIF} were informally introduced in a previous publication~\cite{Rieck17}. Here, we give
a more formal introduction, demonstrate that \glspl{PIF} can be easily and rapidly calculated,
derive several properties that are advantageous for topological data analysis as well as machine
learning, and describe example usage scenarios, such as hypothesis testing and classification.
We make our implementation, experiments, and data publicly
available\footnote{\url{https://github.com/Submanifold/topological-machine-learning}}.

\section{Related Work}

The \emph{persistence curve} is a precursor to \glspl{PIF} that is widely used
in the analysis of Morse--Smale complexes~\cite{Bremer04, Guenther14, Laney06}.
It counts the number of certain critical points, such as minima or maxima, that either occur at
a given persistence threshold or at given point in time.
The curve is then used to determine a relevant threshold, or cut-off parameter for the
simplification of the critical points of a function. To the best of our knowledge, no standadized
variant of these curves appears to exist.

Recognizing that persistence diagrams can be analyzed at multiple scales as well in order to
facilitate hierarchical comparisons, there are some approaches that provide approximations to
persistence diagrams based on, e.g., a smoothing parameter. Among these, the stable kernel of
Reininghaus et al.~\cite{Reininghaus15} is particularly suited for topological machine learning.
Another approach by Adams et al.~\cite{Adams17} transforms a persistence diagram into
a finite-dimensional vector by means of a probability distribution.
Both methods require choosing a set of parameters~(for kernel computations), while
\glspl{PIF} are fully parameter-free. Moreover, \glspl{PIF} also permit other applications, such as
mean calculations and statistical hypothesis testing, which pure kernel methods cannot provide.

Recently, Bubenik~\cite{Bubenik15} introduced \emph{persistence landscapes}, a functional summary of
persistence diagrams. Within his framework, \glspl{PIF} can be considered to represent a summary~(or
projection) of the \emph{rank function}. Our definition of \glspl{PIF} is more straightforward and
easier to implement, however.
Since \glspl{PIF} share several properties of persistence
landscapes---most importantly the existence of simple function-space distance measures---this paper
uses similar experimental setups as Bubenik~\cite{Bubenik15} and Chazal et al.~\cite{Chazal13}.

\section{Persistence Indicator Functions~(PIFs)}

\begin{figure}[tbp]
  \centering
  \subfloat[\label{sfig:Persistence diagram example}]{%
    \iffinal
      \includegraphics{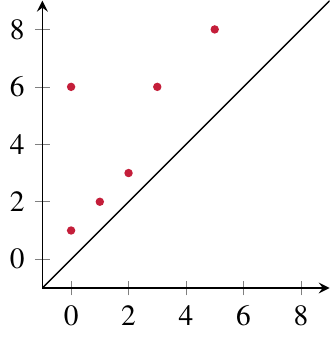}
    \else
      \begin{tikzpicture}
        \begin{axis}[%
          axis x line = bottom,
          axis y line = left,
          xmin        = -1.0,
          xmax        =  9.0,
          ymin        = -1.0,
          ymax        =  9.0,
          width       =  4.5cm,
          height      =  4.5cm,
          mark size   = 1pt,
        ]
          \addplot[cardinal, only marks] file {Data/Illustrations/Persistence_diagram.txt};
          \addplot[domain={-1.0:9.0}, no marks] {x};
        \end{axis}
      \end{tikzpicture}
    \fi
  }
  \subfloat[\label{sfig:Persistence barcode example}]{%
    \iffinal
      \includegraphics{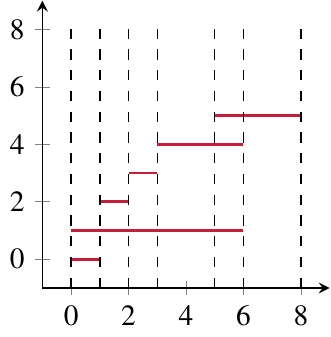}
    \else
      \begin{tikzpicture}
        \begin{axis}[%
          axis x line = bottom,
          axis y line = left,
          xmin        = -1.0,
          xmax        =  9.0,
          ymin        = -1.0,
          ymax        =  9.0,
          width       =  4.5cm,
          height      =  4.5cm,
          mark size   = 1pt,
        ]
          \addplot[thick, cardinal] gnuplot [raw gnuplot] {%
            plot "Data/Illustrations/Persistence_barcode.txt" index 0 with lines;
        };
          \addplot[dashed] gnuplot [raw gnuplot] {%
            plot "Data/Illustrations/Persistence_barcode.txt" index 1 with lines;
        };
        \end{axis}
      \end{tikzpicture}
    \fi
  }
  \subfloat[\label{sfig:Persistence indicator function example}]{%
    \iffinal
      \includegraphics{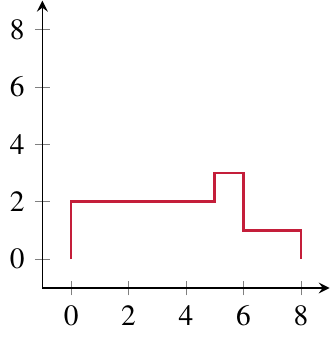}
    \else
      \begin{tikzpicture}
        \begin{axis}[%
          axis x line = bottom,
          axis y line = left,
          xmin        = -1.0,
          xmax        =  9.0,
          ymin        = -1.0,
          ymax        =  9.0,
          width       =  4.5cm,
          height      =  4.5cm,
          mark size   = 1pt,
        ]
        \addplot[thick, cardinal] gnuplot [raw gnuplot]{%
          plot "Data/Illustrations/Persistence_indicator_function.txt" with steps;
        };
        \end{axis}
      \end{tikzpicture}
    \fi
  }
  \caption{%
    A persistence diagram~\protect\subref{sfig:Persistence diagram example}, its persistence
    barcode~\protect\subref{sfig:Persistence barcode example}, and its corresponding persistence
    indicator function~\protect\subref{sfig:Persistence indicator function example}.
    Please note that the interpretation of the axes changes for each plot.
  }
  \label{fig:Persistence indicator function}
\end{figure}

Given a persistence diagram~$\diagram$, i.e., a descriptor of the topological activity of a data
set~\cite{Edelsbrunner10}, we can summarize topological features by calculating an associated
\acrlong{PIF} of $\diagram$ as
\begin{equation}
  \begin{aligned}
    \pif\colon \real &\longrightarrow \mathds{N}\\
    \epsilon &\longmapsto \big| \big\{ (c,d) \in \diagram \mid  \epsilon \in [c,d] \big\} \big|,
  \end{aligned}
\end{equation}
i.e., the number of points in the persistence diagram that, when being treated as a closed interval,
contain the given parameter~$\epsilon$. Equivalently, a \gls{PIF} can be considered to describe
the rank of the $p$\th homology group of a filtration of a simplicial complex.
A \gls{PIF} thus measures the amount of topological activity as a function of the threshold
parameter~$\epsilon$. This parameter is commonly treated as the ``range'' of a function defined
\emph{on} the given data set, e.g., a distance function~\cite{Edelsbrunner02} or an elevation
function~\cite{Agarwal06}. Figure~\ref{fig:Persistence indicator function} demonstrates how to
calculate the persistence indicator function~$\pif(\cdot)$ from a persistence diagram or, equivalently, from
a persistence barcode. For the latter, the calculation becomes particularly easy. In the barcode,
one only has to check the number of intervals that are intersected at any given time by
a vertical line for some value of~$\epsilon$.

\subsection{Properties}

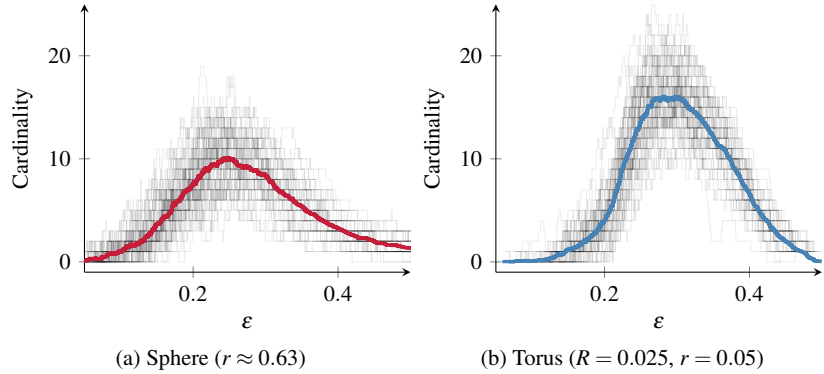
\begin{figure}[tbp]
  \centering
  \subfloat[Sphere~($r \approx 0.63$)]{%
    \iffinal
      \includegraphics{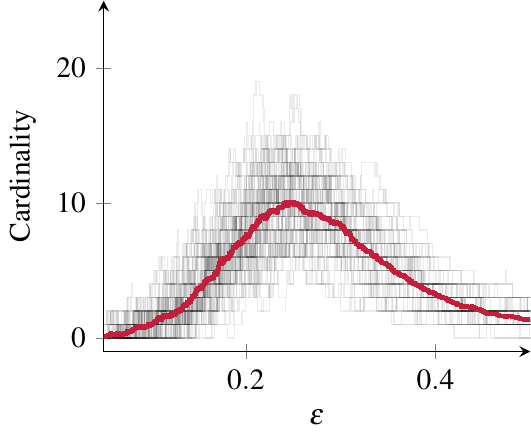}
    \else
      \begin{tikzpicture}
        \begin{axis}[
          axis x line = bottom,
          axis y line = left,
          xmin= 0.05, xmax=00.50,
          ymin=-1.00, ymax=25.00,
          enlargelimits = false,
          scale         = 0.80,
          axis x line   = bottom,
          axis y line   = left,
          xlabel        = $\epsilon$,
          ylabel        = {Cardinality},
          width = 7cm,
        ]
          \pgfplotsinvokeforeach{1,...,50}{%
            \addplot[opacity=0.075] gnuplot[raw gnuplot]{%
              input = sprintf("Data/PIF_sphere_\%02d.txt", #1);
              plot input with steps;
            };
          }
          \addplot[cardinal, very thick] gnuplot[raw gnuplot]{%
            plot "Data/PIF_sphere_mean.txt" with steps;
          };
        \end{axis}
      \end{tikzpicture}
    \fi
  }
  \subfloat[Torus~($R = 0.025$, $r = 0.05$)]{%
    \iffinal
      \includegraphics{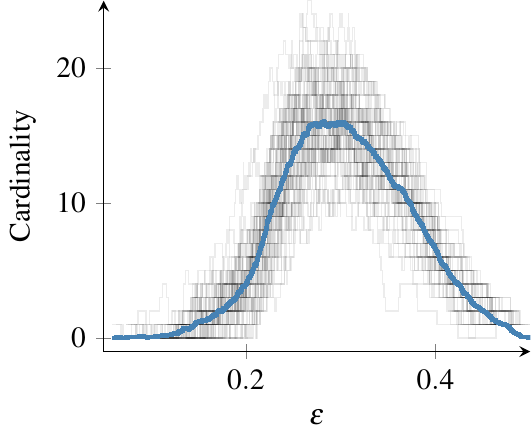}
    \else
      \begin{tikzpicture}
        \begin{axis}[
          axis x line = bottom,
          axis y line = left,
          xmin= 0.05, xmax=00.50,
          ymin=-1.00, ymax=25.00,
          enlargelimits = false,
          scale         = 0.80,
          axis x line   = bottom,
          axis y line   = left,
          xlabel        = $\epsilon$,
          ylabel        = {Cardinality},
          width = 7cm,
        ]
          \pgfplotsinvokeforeach{1,...,50}{%
            \addplot[opacity=0.075] gnuplot[raw gnuplot]{%
              input = sprintf("Data/PIF_torus_\%02d.txt", #1);
              plot input with steps;
            };
          }
          \addplot[yale, very thick] gnuplot[raw gnuplot]{%
            plot "Data/PIF_torus_mean.txt" with steps;
          };
        \end{axis}
      \end{tikzpicture}
    \fi
  }
  \caption{%
    Mean persistence indicator function of the one-dimensional persistence diagrams of a sphere and
    of a torus. Both data sets have been sampled at random and are scaled such that their volume is
    the same.
  }
  \label{fig:Persistence indicator function mean}
\end{figure}

We first observe that the \gls{PIF} only changes at finitely many points. These are given by the
creation and destruction times, i.e., the $x$- and $y$-coordinates, of points in the persistence
diagram. The \gls{PIF} may thus be written as a sum of appropriately scaled \emph{indicator
functions}~(hence the name) of the form~$\pif[\mathcal{I}](\cdot)$ for some interval~$\mathcal{I}$. Within
the interval~$\mathcal{I}$, the value of $\pif(\cdot)$ does \emph{not} change. Hence, the \gls{PIF} is a \emph{step function}.
Since step functions are compatible with addition and scalar multiplication, \glspl{PIF} form
a vector space. The addition of two \glspl{PIF} corresponds to calculating the union of their
corresponding persistence diagrams, while taking multiplicities of points into account. As
a consequence, we can calculate the \emph{mean} of a set of \glspl{PIF} $\{\pif^1,\dots,\pif^n\}$ as
\begin{equation}
  \overline{\pif}(\cdot) := \frac{1}{n} \sum_{i=1}^{n} \pif^i(\cdot),
\end{equation}
i.e., the standard pointwise mean of set of elements. In contrast to persistence diagrams, for which
a mean is not uniquely defined and hard to calculate~\cite{Turner14}, this calculation only involves
addition~(union) and scalar multiplications of sets of intervals.
Figure~\ref{fig:Persistence indicator function mean} depicts mean persistence indicator functions
for two randomly-sampled data sets. We see that the resulting mean persistence indicator functions
already introduce a visual separation between the two data sets.

As a second derived property, we note that the absolute value of a step function~(and that of
a \gls{PIF}) always exists; one just calculates the absolute value for every interval in which the
step function does \emph{not} change. The absolute value of a step function is again a step
function, so the Riemann integral of a \gls{PIF} is well-defined, giving rise to their $1$-norm
as
\begin{equation}
  \| \pif \|_1 := \int_\real|\pif(x)|\d{x},
\end{equation}
which is just the standard norm of an \Lp[1]-space. The preceding equation requires the use of an
absolute value because linear operations on \glspl{PIF} may result in negative values.
The integral of a \gls{PIF}~(or its absolute
value) decomposes into a sum of integrals of individual step functions, defined over some interval
$[a,b]$. Letting the value of the step function over this interval be $l$, the integral of the step
function is given as $l\cdot|b-a|$, i.e., the volume of the interval scaled by the value the step function
assumes on it. We can also extend this norm to that of an \Lp-space, where $p \in \real$.
To do so, we observe that the $p$\th power of any step function is
well-defined---we just raise the value it takes to the $p$\th power. Consequently, the $p$\th power
of a \gls{PIF} is also well-defined and we define
\begin{equation}
  \| \pif \|_p := \Big(\int_\real|\pif(x)|^p\d{x}\Big)^\frac{1}{p},
\end{equation}
which is the standard norm of an \Lp-space. Calculating this integral again involves decomposing the
range of the \gls{PIF} into individual step functions and calculating their integral.
We have $\| \pif \|_p < \infty$ for every $p \in \real$ because there are only finitely many points
in a persistence diagram, so the integrals of the individual step functions involved in the norm
calculation are bounded.

\runinhead{Hypothesis Testing}

Treating the norm of a \gls{PIF} as a random variable, we can perform topology-based hypothesis
testing similar to persistence landscapes~\cite{Bubenik15}.
Given two different samples of persistence diagrams, $\{\diagram^1_1, \dots,
\diagram^1_n\}$ and $\{\diagram^2_1, \dots, \diagram^2_n\}$, we calculate the $1$-norm of their
respective mean \glspl{PIF} as $Y_1$ and $Y_2$, and the variances
\begin{equation}
  \sigma_i^2 := \frac{1}{n-1} \sum_{j=1}^{n}\Big(\big|\pif[\diagram^i_j]\big| - Y_i\Big)^2,
\end{equation}
for $i \in \{1,2\}$. We may then perform a standard two-sample $z$-test to check whether the two
means are likely to be the same. To this end, we calculate the $z$-score as
\begin{equation}
  z := \frac{Y_1 - Y_2}{\sqrt{s_1^2/n-s_2^2/n}}
\end{equation}
and determine the critical values at the desired $\alpha$-level, i.e., the significance level, from
the quantiles of a normal distribution. If $z$ is outside the interval spanned by the critical
values, we \emph{reject} the null hypothesis, i.e., we consider the two means to be different.
For the example shown in Figure~\ref{fig:Persistence indicator function mean}, we obtain $Y_1
\approx 2.135$, $s_1^2 \approx 0.074$, $Y_2 \approx 2.79$, and $s_2^2 \approx 0.093$.
Using $\alpha = 0.01$, the critical values are given by $z_1 \approx -2.58$ and $z_2 \approx
2.58$. Since $z \approx -11.09$, we \emph{reject} the null hypothesis with $p \approx 1.44 \times
10^{-28} \ll 0.01$. Hence, \glspl{PIF} can be used to confidently discern random samples of
a sphere from those of a torus with the same volume.

\runinhead{Stability}
%
The $1$-norm of a \gls{PIF} is connected to the \emph{total persistence}~\cite{Cohen-Steiner10},
i.e., the sum of all persistence values in a persistence diagram. We have
\begin{equation}
  \| \pif \|_1 = \int_\real|\pif(x)|\d{x} = \sum_{I\in\mathcal{I}} c_I \volume(I),
\end{equation}
where $\mathcal{I}$ denotes a set of intervals for which the number of active persistence pairs does
not change, and $c_I$ denotes their count.
We may calculate this partition from a persistence barcode, as shown in Figure~\ref{sfig:Persistence
barcode example}, by going over all the dotted slices, i.e., the intervals between pairs of
start and endpoints of each interval.
The sum over all these volumes is equal to the total persistence of the set of intervals, because we
can split up the volume calculation of a single interval over many sub-interval and, in total, the
volume of every interval is only accumulated once. Hence,
\begin{equation}
  \| \pif \|_1  \int_\real|\pif(x)|\d{x} = \sum_{(c,d)\in\diagram} |d-c| = \sum_{(c,d)\in\diagram}\pers(c,d) = \pers(\diagram),
\end{equation}
where $\pers(\diagram)$ denotes the total persistence of the persistence diagram. According
to a stability theorem by Cohen-Steiner et al.~\cite{Cohen-Steiner10}, the $1$-norm of a \gls{PIF}
is thus stable with respect to small perturbations.
We leave the derivation of a similar equation for the general \Lp-norm of a \gls{PIF} for future work.

\subsection{The Bootstrap for Persistence Indicator Functions}

Developed by Efron and Tibshirani~\cite{Efron93}, the bootstrap is a general statistical method
for---among other applications---computing confidence intervals. We give a quick and cursory
introduction before showing how this method applies to persistence indicator functions;
please refer to Chazal et al.~\cite{Chazal13} for more details.

Assume that we have a set of independent and identically distributed variables $X_1$, \dots,
$X_n$, and we want to estimate a real-valued parameter $\theta$ that corresponds to their
distribution.
Typically, we may estimate $\theta$ using a statistic $\widehat{\theta} := s(X_1,\dots,X_n)$, i.e.,
some function of the data. A common example is to use $\theta$ as the population mean, while
$\widehat{\theta}$ is the sample mean. If we want to calculate \emph{confidence intervals} for our
estimate $\widehat{\theta}$, we require the distribution of the difference $\theta
- \widehat{\theta}$. This distribution, however, depends on the unknown distribution of the
variables, so we have to approximate it using an empirical distribution.
Let $X_1^\ast$, \dots, $X_n^\ast$ be a sample of the original variables, drawn
with replacement. We can calculate $\widehat{\theta}^\ast := s(X_1^\ast,\dots,X_n^\ast)$ and
approximate the unknown distribution by the empirical distribution of $\widehat{\theta}
- \widehat{\theta}^\ast$, which, even though it is not computable analytically, can be
approximated by repeating the sampling procedure a sufficient number of times. The quantiles of the
approximated distribution may then be used to construct confidence intervals, leading
to the following method:
\begin{enumerate}
  \item Calculate an estimate of $\theta$ from the input data using $\widehat{\theta} := s(X_1,\dots,X_n)$.
  \item Obtain $X_1^\ast$, \dots, $X_n^\ast$~(sample \emph{with} replacement) and calculate $\widehat{\theta}^\ast := s(X_1^\ast, \dots, X_n^\ast)$.
  \item Repeat the previous step $B$ times to obtain $\widehat{\theta}^\ast_1$, \dots, $\widehat{\theta}^\ast_B$.
  \item Given $\alpha$, compute an approximate $(1-2\alpha)$ quantile interval as
        \begin{equation}
          [\widehat{\theta}_1, \widehat{\theta}_2] \approx [\widehat{\theta}^{\ast(\alpha)}_B, \widehat{\theta}^{\ast(1-\alpha)}_B],
        \end{equation}
        where $\widehat{\theta}^{\ast(\alpha)}_B$ refers to the $\alpha$\th empirical quantile of the bootstrap replicates from the previous step.
\end{enumerate}
This procedure yields both a lower bound and an upper bound for $\widehat{\theta}$. It is also
referred to as the percentile interval for bootstraps~\cite[pp.\ 168--176]{Efron93}. More
complicated versions---yielding ``tighter'' confidence intervals---of this procedure exist, but the
basic method of sampling with replacement and computing the statistic~$s(\cdot)$ on the bootstrap
samples remains the same.

\begin{figure}[tbp]
  \centering
  \tikzset{external/system call = {%
      lualatex \tikzexternalcheckshellescape -halt-on-error -interaction=batchmode -jobname "\image" "\texsource"
    }
  }
  \subfloat[Sphere\label{sfig:Mean persistence indicator function sphere}]{%
    \iffinal
      \includegraphics{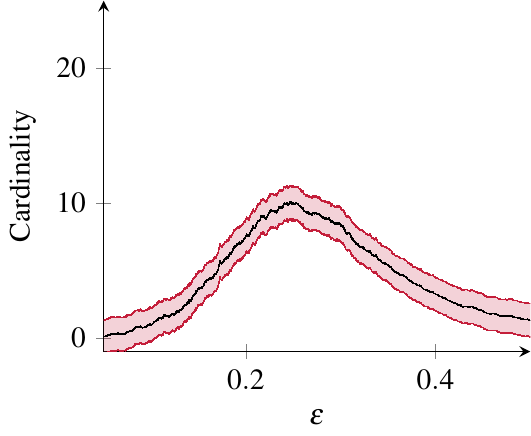}
    \else
      \begin{tikzpicture}
        \centering
        \begin{axis}[
          xmin= 0.05, xmax=00.50,
          ymin=-1.00, ymax=25.00,
          enlargelimits = false,
          scale         = 0.80,
          axis x line   = bottom,
          axis y line   = left,
          xlabel        = $\epsilon$,
          ylabel        = {Cardinality},
          width = 7cm,
        ]
          \addplot[name path=F, black] gnuplot[raw gnuplot]{%
              plot "Data/PIF_sphere_mean_plus_confidence.txt" index 0 with steps;
          };
          \addplot[name path=U, cardinal] gnuplot[raw gnuplot]{%
              plot "Data/PIF_sphere_mean_plus_confidence.txt" index 1 with steps;
          };
          \addplot[name path=L, cardinal] gnuplot[raw gnuplot]{%
              plot "Data/PIF_sphere_mean_plus_confidence.txt" index 2 with steps;
          };
          \addplot[cardinal, fill opacity=0.2] fill between [of=U and L];
        \end{axis}
      \end{tikzpicture}
    \fi
  }
  \subfloat[Torus\label{sfig:Mean persistence indicator function torus}]{%
    \iffinal
      \includegraphics{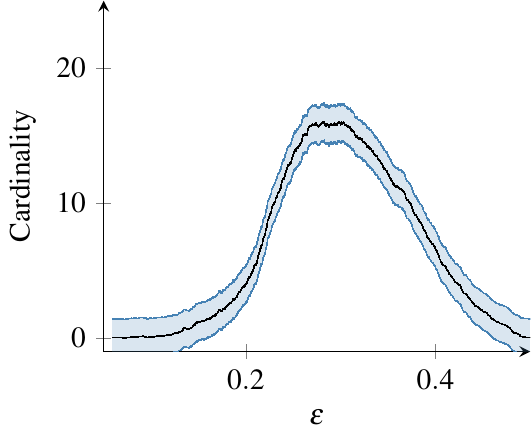}
    \else
      \begin{tikzpicture}
        \centering
        \begin{axis}[
          xmin= 0.05, xmax=00.50,
          ymin=-1.00, ymax=25.00,
          enlargelimits = false,
          scale         = 0.80,
          axis x line   = bottom,
          axis y line   = left,
          xlabel        = $\epsilon$,
          ylabel        = {Cardinality},
          width = 7cm,
        ]
          \addplot[name path=F, black] gnuplot[raw gnuplot]{%
              plot "Data/PIF_torus_mean_plus_confidence.txt" index 0 with steps;
          };
          \addplot[name path=U, yale] gnuplot[raw gnuplot]{%
              plot "Data/PIF_torus_mean_plus_confidence.txt" index 1 with steps;
          };
          \addplot[name path=L, yale] gnuplot[raw gnuplot]{%
              plot "Data/PIF_torus_mean_plus_confidence.txt" index 2 with steps;
          };
          \addplot[yale, fill opacity=0.2] fill between [of=U and L];
        \end{axis}
      \end{tikzpicture}
    \fi
  }
  \caption{
    Confidence bands at the $\alpha = 0.05$ level for the mean persistence indicator functions of
    a sphere and a torus. The confidence band is somewhat tighter for $\epsilon \geq 0.2$.
  }
  \label{fig:Confidence bands sphere torus}
\end{figure}
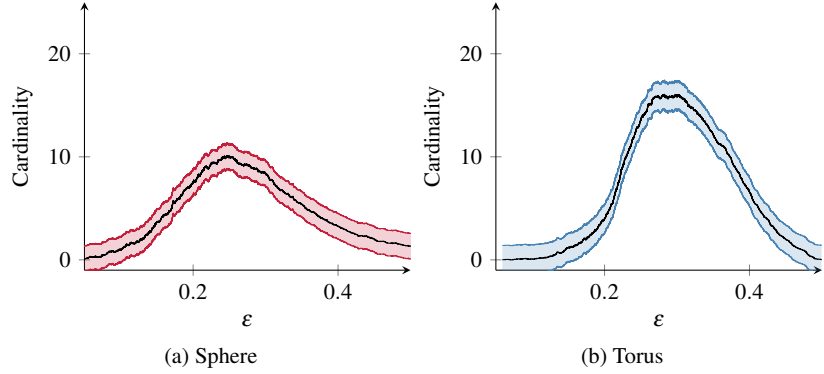

In order to apply the bootstrap to \emph{functions}, we require empirical processes~\cite{Kosorok08}.
The goal is to find a \emph{confidence band} for a function~$f(x)$, i.e., a pair of functions
$l(x)$ and $u(x)$ such that the probability that $f(x) \in [l(x), u(x)]$ for $x \in \real$ is at
least $1-\alpha$. Given a function $f$, let $Pf := \int f \d{P}$ and $P_n f := n^{-1} \sum_{i=1}^{n}
f(X_i)$. We obtain a \emph{bootstrap empirical process} as
\begin{equation}
  \{ \empiricalsample f \}_{f\in\mathcal{F}} := \{ \sqrt{n} \left(P_n^\ast f -  P_n f \right) \},
\end{equation}
where $P_n^\ast := n^{-1} \sum_{i=1}^{n} f(X_i^\ast)$ is defined on the bootstrap samples~(as
introduced above). Given the convergence of this empirical process, we may calculate
\begin{equation}
  \widehat{\theta} := \sup_{f\in\mathcal{F}} | \empiricalsample f |,
\end{equation}
which yields a statistic to use for the bootstrap as defined above. From the corresponding
quantile, we ultimately obtain $[\widehat{\theta}_1, \widehat{\theta}_2]$ and calculate a confidence band
\begin{equation}
  C_n(f) := \left[ P_n f - \frac{\widehat{\theta}_1}{n}, P_n f + \frac{\widehat{\theta}_2}{n} \right]
\end{equation}
for the empirical mean of a set of \glspl{PIF}. Figure~\ref{fig:Confidence bands sphere torus}
depicts an example of confidence bands for the mean \gls{PIF} of the sphere and torus samples. We can
see that the confidence band for the torus is extremely tight for $\epsilon \in [0.2,0.3]$,
indicating that the limit behavior of samples from a torus is different at this scale from the limit
behavior of samples from a sphere.

\subsection{Distances and Kernels}
\label{sec:Distances}

Given two persistence diagrams $\diagram_i$ and $\diagram_j$, we are often interested in their
dissimilarity~(or distance). Having seen that linear combinations and norms of \glspl{PIF} are
well-defined, we can define a family of distances as
\begin{equation}
  \dist_p(\pif[\diagram_i], \pif[\diagram_j]) := \Big(\int_\real|\pif[\diagram_i](x)-\pif[\diagram_j](x)|^p\d{x}\Big)^\frac{1}{p},
\end{equation}
with $p\in\real$. Since the norm of a \gls{PIF} is well-defined, this expression is a metric
in the mathematical sense. Note that its calculation requires essentially only evaluating all
individual step functions of the difference of the two \glspl{PIF} once. Hence, its complexity is
\emph{linear} in the number of sub-intervals.

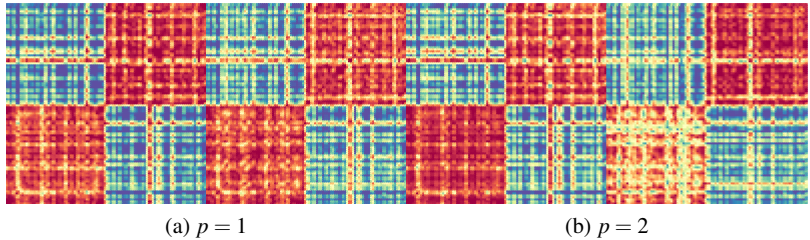
\begin{figure}[tbp]
  \centering
  \tikzset{external/system call = {%
      lualatex \tikzexternalcheckshellescape -halt-on-error -interaction=batchmode -jobname "\image" "\texsource"
    }
  }
  \subfloat[$p=1$]{%
    \iffinal
      \includegraphics{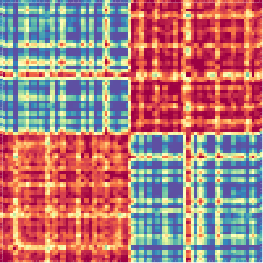}%
      \includegraphics{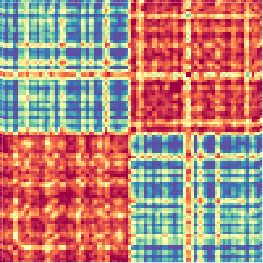}%
    \else
      \begin{tikzpicture}
        \begin{axis}[%
          enlargelimits     = true,
          axis x line       = none,
          axis y line       = none,
          scale only axis   = true,
          unit vector ratio = 1 1 1,
          scale             = 0.75,
        ]
          \addplot[%
            matrix plot*,
            point meta     = explicit,
            mesh/ordering  = colwise,
            faceted color  = gray,
            shader         = flat,
            line width     = 0.01pt,
            point meta min = 0.90, 
            point meta max = 1.75, 
            colormap/Spectral,
            ]
            table[header=false, meta index=2] {Data/PIF1_distances_sphere_torus_transformed.txt};
        \end{axis}
      \end{tikzpicture}
      \begin{tikzpicture}
        \begin{axis}[%
          enlargelimits     = true,
          axis x line       = none,
          axis y line       = none,
          scale only axis   = true,
          unit vector ratio = 1 1 1,
          scale             = 0.75,
        ]
          \addplot[%
            matrix plot*,
            point meta     = explicit,
            mesh/ordering  = colwise,
            faceted color  = gray,
            shader         = flat,
            line width     = 0.01pt,
            point meta min = 0.80, 
            point meta max = 1.20, 
            colormap/Spectral,
            ]
            table[header=false, meta index=2] {Data/W1_distances_sphere_torus_transformed.txt};
        \end{axis}
      \end{tikzpicture}
    \fi
  }
  \subfloat[$p=2$]{%
    \iffinal
      \includegraphics{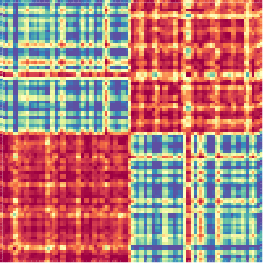}%
      \includegraphics{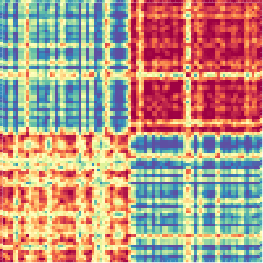}%
    \else
      \begin{tikzpicture}
        \begin{axis}[%
          enlargelimits     = true,
          axis x line       = none,
          axis y line       = none,
          scale only axis   = true,
          unit vector ratio = 1 1 1,
          scale             = 0.75,
        ]
          \addplot[%
            matrix plot*,
            point meta     = explicit,
            mesh/ordering  = colwise,
            faceted color  = gray,
            shader         = flat,
            line width     = 0.01pt,
            point meta min = 3,  
            point meta max = 12, 
            colormap/Spectral,
            ]
            table[header=false, meta index=2] {Data/PIF2_distances_sphere_torus_transformed.txt};
        \end{axis}
      \end{tikzpicture}
      \begin{tikzpicture}
        \begin{axis}[%
          enlargelimits     = true,
          axis x line       = none,
          axis y line       = none,
          scale only axis   = true,
          unit vector ratio = 1 1 1,
          scale             = 0.75,
          colormap/Spectral,
        ]
          \addplot[%
            matrix plot*,
            point meta    = explicit,
            mesh/ordering = colwise,
            faceted color = gray,
            shader        = flat,
            line width    = 0.01pt,
            point meta min = 0.38, 
            point meta max = 0.47, 
            ]
            table[header=false, meta index=2] {Data/W2_distances_sphere_torus_transformed.txt};
        \end{axis}
      \end{tikzpicture}
    \fi
  }
  \caption{%
    A comparison of distance matrices obtained using the distance measure~$\dist_p$ for \glspl{PIF},
    and the corresponding $p$\th Wasserstein distance~$\wasserstein_p$. The left matrix of each
    group shows~$\dist_p$, while the right matrix depicts~$\wasserstein_p$.
    Red indicates close objects~(small distances), while blue indicates far objects~(large distances).
  }
  \label{fig:Distance comparison sphere torus}
\end{figure}

\runinhead{Example}
%
Figure~\ref{fig:Distance comparison sphere torus} depicts pairwise distance matrices for random
samples of a sphere and of a torus. The first matrix of each group is obtained via $\dist_p$ for
\glspl{PIF}, while the second matrix in each group is obtained by the $p$\th Wasserstein distance.
We observe two groups of data sets in all matrices, meaning that both classes of distances
are suitable for detecting differences.\\

\noindent We can also employ the distance defined above to obtain
a \emph{kernel}~\cite{Schoelkopf02}. To this end, we define
\begin{equation}
  k_p(\diagram_i, \diagram_j) := -\dist_p(\pif[\diagram_i], \pif[\diagram_j]),
\end{equation}
where $p \in \{1,2\}$ because we need to make sure that the kernel is \emph{conditionally positive
definite}~\cite{Schoelkopf02}. This kernel permits using \glspl{PIF} with many modern machine
learning algorithms.
As an illustrative example, we will use \emph{kernel support vector machines}~\cite{Schoelkopf02}
to separate random samples of a sphere and a torus.

\runinhead{Example}
%
Again, we use 100 random samples~(50 per class) from a sphere and a torus. We only use
one-dimensional persistence diagrams, from which we calculate \glspl{PIF}, from which we then obtain
pairwise kernel matrices using both $k_1$ and $k_2$. Finally, we train a support vector machine
using nested stratified $5$-fold cross-validation.
With $k_1$, we obtain an average accuracy of $0.98 \pm 0.049$, whereas with $k_2$, we obtain an
average accuracy of $0.95 \pm 0.063$. The decrease in accuracy is caused by the additional smoothing
introduced in this kernel.

\section{Applications}

In the following, we briefly discuss some potential application scenarios for \glspl{PIF}. We use
only data sets that are openly available in order to make our results comparable. For all
machine learning methods, we use \textsc{scikit-learn}~\cite{Pedregosa11}.

\subsection{Analysis of Random Complexes}

\begin{figure}[tbp]
  \centering
  \tikzset{external/system call = {%
      lualatex \tikzexternalcheckshellescape -halt-on-error -interaction=batchmode -jobname "\image" "\texsource"
    }
  }
  \iffinal
    \includegraphics{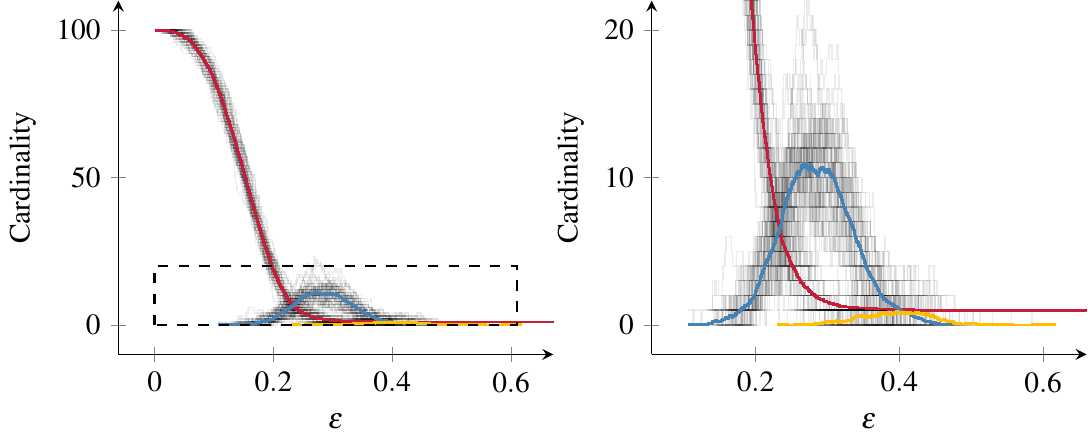}
  \else
    \begin{tikzpicture}
      \begin{groupplot}[%
        group style={group size=2 by 1},
        axis x line   = bottom,
        axis y line   = left,
        enlargelimits = true,
        xmax          = 0.61,
        width         = 6cm,
        xlabel        = $\epsilon$,
        ylabel        = {Cardinality},
      ]
        \nextgroupplot
          \pgfplotsinvokeforeach{0,...,49}{%
            \addplot[opacity=0.075] gnuplot[raw gnuplot] {%
              input = sprintf("Data/Random_complexes/PIF_\%02d_d0.txt", #1);
              plot input with steps;
            };
          }
          \pgfplotsinvokeforeach{0,...,49}{%
            \addplot[opacity=0.075] gnuplot[raw gnuplot] {%
              input = sprintf("Data/Random_complexes/PIF_\%02d_d1.txt", #1);
              plot input with steps;
            };
          }
          \pgfplotsinvokeforeach{0,...,49}{%
            \addplot[opacity=0.075] gnuplot[raw gnuplot] {%
              input = sprintf("Data/Random_complexes/PIF_\%02d_d2.txt", #1);
              plot input with steps;
            };
          }
          \addplot[cardinal,thick] gnuplot[raw gnuplot] {%
            plot "Data/Random_complexes/PIF_mean_d0.txt" with steps;
          };
          \addplot[yale,thick] gnuplot[raw gnuplot] {%
            plot "Data/Random_complexes/PIF_mean_d1.txt" with steps;
          };
          \addplot[amber,thick] gnuplot[raw gnuplot] {%
            plot "Data/Random_complexes/PIF_mean_d2.txt" with steps;
          };
          \draw [black,dashed,thick] (0,0) rectangle (0.61,20);
        \nextgroupplot[ymax = 20]
          \pgfplotsinvokeforeach{0,...,49}{%
            \addplot[opacity=0.075] gnuplot[raw gnuplot] {%
              input = sprintf("Data/Random_complexes/PIF_\%02d_d0.txt", #1);
              plot input with steps;
            };
          }
          \pgfplotsinvokeforeach{0,...,49}{%
            \addplot[opacity=0.075] gnuplot[raw gnuplot] {%
              input = sprintf("Data/Random_complexes/PIF_\%02d_d1.txt", #1);
              plot input with steps;
            };
          }
          \pgfplotsinvokeforeach{0,...,49}{%
            \addplot[opacity=0.075] gnuplot[raw gnuplot] {%
              input = sprintf("Data/Random_complexes/PIF_\%02d_d2.txt", #1);
              plot input with steps;
            };
          }
          \addplot[cardinal,thick] gnuplot[raw gnuplot] {%
            plot "Data/Random_complexes/PIF_mean_d0.txt" with steps;
          };
          \addplot[yale,thick] gnuplot[raw gnuplot] {%
            plot "Data/Random_complexes/PIF_mean_d1.txt" with steps;
          };
          \addplot[amber,thick] gnuplot[raw gnuplot] {%
            plot "Data/Random_complexes/PIF_mean_d2.txt" with steps;
          };
      \end{groupplot}
    \end{tikzpicture}
  \fi
  \caption{%
    \glspl{PIF} for random complexes sampled over a unit cube in $\real^3$. Dimensions zero~(red),
    one~(blue), and two~(yellow) are shown. To show the peak in dimension two better, the right-hand
    side shows a ``zoomed'' version of the first chart~(dashed region).
  }
  \label{fig:Random box}
\end{figure}
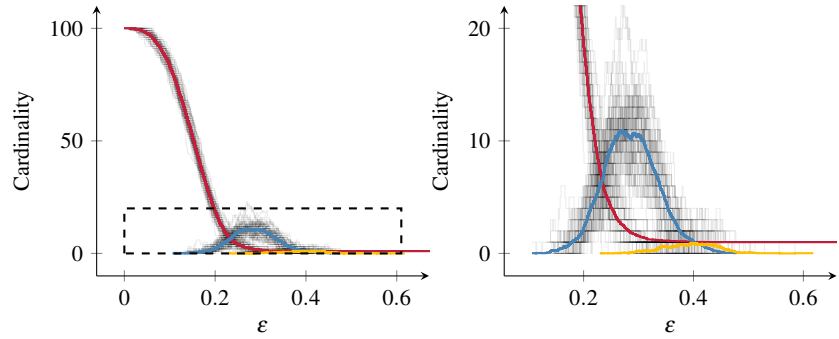

It is often useful to know to what extent a data set exhibits random fluctuations. To this end, we
sampled 100 points from a unit cube in $\real^3$, which has the topology of a point, i.e., no
essential topological features in dimensions $> 0$.
We calculated the Vietoris--Rips complex at a scale such that no essential topological features
remain in dimensions $\geq 1$, and obtained \glspl{PIF}, which are depicted in
Figure~\ref{fig:Random box} along with their corresponding mean.
All functions use a common axis in order to simplify their comparison.
We first comment on the dynamics of these data. Topological activity is ``shifted'', meaning that
topological features with a different dimensionality are not active at the same scale.  The maximum
of topological activity in dimension one~(blue curve) is only reached when there are few
zero-dimensional features. The maximum in dimension two~(yellow curve) also does not coincide with
the maximum in dimension one.
These results are consistent with a limit theorem of Bobrowski and Kahle~\cite{Bobrowski14}, who
showed that~(persistent) Betti numbers follow a Poisson distribution.

\begin{figure}[tbp]
  \centering
  \tikzset{external/system call = {%
      lualatex \tikzexternalcheckshellescape -halt-on-error -interaction=batchmode -jobname "\image" "\texsource"
    }
  }
  \iffinal
    \includegraphics{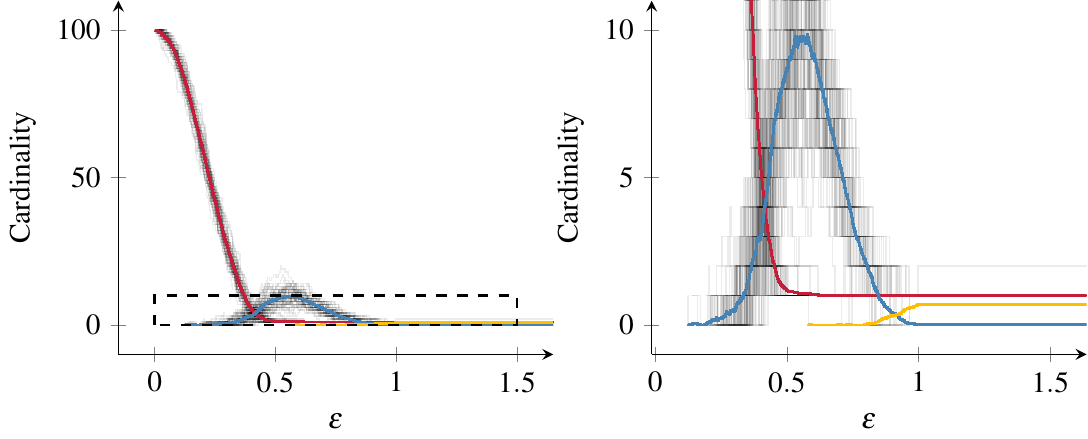}
  \else
    \begin{tikzpicture}
      \begin{groupplot}[%
        group style={group size=2 by 1},
        axis x line   = bottom,
        axis y line   = left,
        enlargelimits = true,
        xmax          = 1.5,
        width         = 6cm,
        xlabel        = $\epsilon$,
        ylabel        = {Cardinality},
      ]
      \nextgroupplot
        \pgfplotsinvokeforeach{0,...,49}{%
          \addplot[opacity=0.075] gnuplot[raw gnuplot] {%
            set xrange [0:1.5];
            set yrange [0:100];
            input = sprintf("Data/Sphere_r10/PIF_\%02d_d0.txt", #1);
            plot input with steps;
          };
        }
        \pgfplotsinvokeforeach{0,...,49}{%
          \addplot[opacity=0.075] gnuplot[raw gnuplot] {%
            set xrange [0:1.5];
            set yrange [0:100];
            input = sprintf("Data/Sphere_r10/PIF_\%02d_d1.txt", #1);
            plot input with steps;
          };
        }
        \pgfplotsinvokeforeach{0,...,49}{%
          \addplot[opacity=0.075] gnuplot[raw gnuplot] {%
            set xrange [0:1.5];
            set yrange [0:100];
            input = sprintf("Data/Sphere_r10/PIF_\%02d_d2.txt", #1);
            plot input with steps;
          };
        }
        \addplot[cardinal, thick] gnuplot[raw gnuplot] {%
          plot "Data/Sphere_r10/PIF_mean_d0.txt" with steps;
        };
        \addplot[yale, thick] gnuplot[raw gnuplot] {%
          plot "Data/Sphere_r10/PIF_mean_d1.txt" with steps;
        };
        \addplot[amber, thick] gnuplot[raw gnuplot] {%
          plot "Data/Sphere_r10/PIF_mean_d2.txt" with steps;
        };
        \draw [black,dashed,thick] (0,0) rectangle (1.5,10);
      \nextgroupplot[ymax = 10]
        \pgfplotsinvokeforeach{0,...,49}{%
          \addplot[opacity=0.075] gnuplot[raw gnuplot] {%
            set xrange [0:1.5];
            set yrange [0:100];
            input = sprintf("Data/Sphere_r10/PIF_\%02d_d0.txt", #1);
            plot input with steps;
          };
        }
        \pgfplotsinvokeforeach{0,...,49}{%
          \addplot[opacity=0.075] gnuplot[raw gnuplot] {%
            set xrange [0:1.5];
            set yrange [0:100];
            input = sprintf("Data/Sphere_r10/PIF_\%02d_d1.txt", #1);
            plot input with steps;
          };
        }
        \pgfplotsinvokeforeach{0,...,49}{%
          \addplot[opacity=0.075] gnuplot[raw gnuplot] {%
            set xrange [0:1.5];
            set yrange [0:100];
            input = sprintf("Data/Sphere_r10/PIF_\%02d_d2.txt", #1);
            plot input with steps;
          };
        }
        \addplot[cardinal, thick] gnuplot[raw gnuplot] {%
          plot "Data/Sphere_r10/PIF_mean_d0.txt" with steps;
        };
        \addplot[yale, thick] gnuplot[raw gnuplot] {%
          plot "Data/Sphere_r10/PIF_mean_d1.txt" with steps;
        };
        \addplot[amber, thick] gnuplot[raw gnuplot] {%
          plot "Data/Sphere_r10/PIF_mean_d2.txt" with steps;
        };
      \end{groupplot}
    \end{tikzpicture}
  \fi
  \caption{%
    \glspl{PIF} for random samples of a sphere with $r = 1.0$. Again, dimensions zero~(red),
    one~(blue), and two~(yellow) are shown, along with a ``zoomed'' version of the first
    chart~(dashed region).
  }
  \label{fig:Random sphere}
\end{figure}
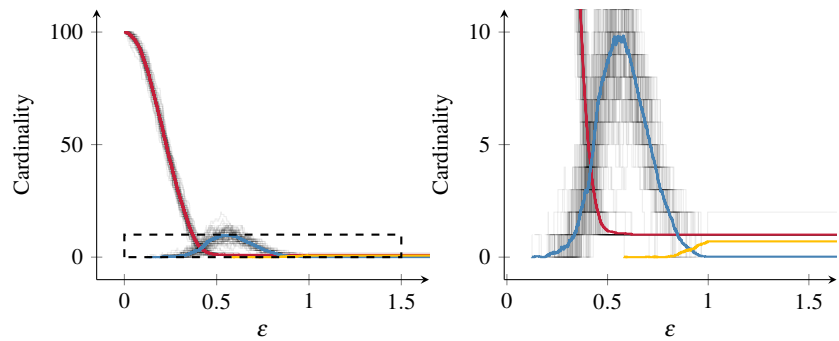

By contrast, for a data set with a well-defined topological structure, such as a $2$-sphere, the
\glspl{PIF} exhibit a different behavior.
Figure~\ref{fig:Random sphere} depicts all \glspl{PIF} of random samples from a sphere. We did not
calculate the Vietoris--Rips complex for all possible values of $\epsilon$ but rather selected an
$\epsilon$ that is sufficiently large to capture the correct Betti numbers of the sphere.
Here, we observe that the stabilization of topological activity in dimension zero roughly coincides
with the maximum of topological activity in dimension one. Topological activity in dimension two
only starts to increase for larger scales, staying stable for a long interval. This activity
corresponds to the two-dimensional void of the $2$-sphere that we detect using persistent homology.

\glspl{PIF} can thus be used to perform a test for ``topological randomness'' in real-world data.
This is useful for deciding whether a topological approximation is suitable or needs
to be changed~(e.g., by calculating a different triangulation, using $\alpha$-shapes, etc.).
Moreover, we can use a \gls{PIF} to detect the presence or absence of a shared ``scale'' in data
sets.
For the random complexes, there is \emph{no} value for $\epsilon$ in which stable topological
activity occurs in more than one dimension, whereas for the sphere, we observe a stabilization
starting from $\epsilon \approx 0.75$.

\subsection{Shakespearean Co-Occurrence Networks}
%
In a previous work, co-occurrence networks from a machine-readable corpus of Shakespeare's
plays~\cite{Rieck16b} have been extracted. Their topological structure under various aspects has
been analyzed, using, for example, their clique communities~\cite{Rieck17} to
calculate two-dimensional embeddings of the individual networks.
The authors observed that comedies form clusters in these embeddings, which indicates that they are
more similar to each other than to plays of another category.
Here, we want to show that it is possible to differentiate between comedies and non-comedies by
using the kernel induced by \glspl{PIF}. Among the 37 available plays, 17 are comedies, giving us
a baseline probability of $0.46$ if we merely ``guess'' the class label of a play. Since the number
of plays is not sufficiently large to warrant a split into test and training data, we use various
cross-validation techniques, such as \emph{leave-one-out}. The reader is referred to
Kohavi~\cite{Kohavi95} for more details.
Table~\ref{tab:Shakespeare results} reports all results; we observe that $k_1$ outperforms $k_2$.
Since $k_2$ emphasizes small-scale differences, the number of topological features in two networks
that are to be compared should be roughly equal.
This is \emph{not} the case for most of the comedies, though. We stress that these results are only
a demonstration of the capabilities of \glspl{PIF}; the comparatively low accuracy is partially due
to the fact that networks were extracted automatically.
It is interesting to note \emph{which} plays tend to be mislabeled.
For $k_1$, \textsc{All's Well That Ends Well}, \textsc{Cymbeline}, and \textsc{The Winter's Tale}
are mislabeled more than all other plays. This is consistent with research by Shakespeare scholars
who suggest different categorization schemes for these~(and other) \emph{problem plays}.

\begin{table}[tbp]
  \caption{%
    Classifier performance for Shakespearean co-occurrence networks. Classification based on $k_1$
    outperforms the second kernel $k_2$.
  }
  \centering
  \setlength{\tabcolsep}{8pt}
  \begin{tabular}{ccccccc}
    \toprule
    Kernel & {$5$-fold} & {LOO}  & {LPO~($p = 2$)} & {LPO~($p = 3$)} & {LPO~($p = 4$)} & {Split}\\
    \midrule
    $k_1$  & 0.84       & 0.83 & 0.83 & 0.80 & 0.79 & 0.80\\
    $k_2$  & 0.64       & 0.00 & 0.67 & 0.68 & 0.68 & 0.68\\
    \bottomrule
  \end{tabular}
  \label{tab:Shakespeare results}
\end{table}

\subsection{Social Networks}

Yanardag and Vishwanathan~\cite{Yanardag15} crawled the popular social news aggregation website
\url{reddit.com} in order to obtain a set of graphs from online discussions. In each graph, the
nodes correspond to users and an edge signifies that a certain user responded to a another user's
comment. The graphs are partitioned into two classes, one of them representing discussion-based
forums~(in which users typically communicate freely among each other), the other representing
communities based on a question--answer format~(in which users typically only respond to the creator
of a topic). The data set is fully-balanced.

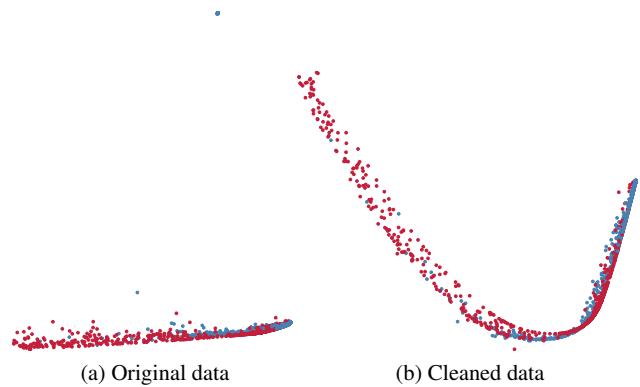
\begin{figure}[tbp]
  \centering
  \subfloat[Original data\label{sfig:reddit PCA original}]{%
    \iffinal
      \includegraphics{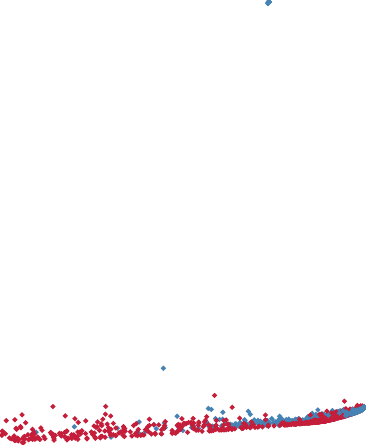}
    \else
      \begin{tikzpicture}
        \begin{axis}[%
          axis lines = none,
          rotate     = -90,
          scatter/classes = {%
            -1={cardinal},
             1={yale}
          },
          height = 6cm,
        ]
          \addplot[scatter, only marks, mark size = 0.5pt, scatter src = explicit symbolic]
            file {Data/REDDIT-BINARY/PCA/Kernel_PCA.txt};
        \end{axis}
      \end{tikzpicture}
    \fi
  }
  \subfloat[Cleaned data\label{sfig:reddit PCA cleaned}]{%
    \iffinal
      \includegraphics{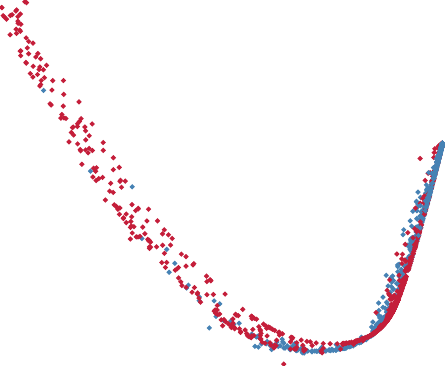}
    \else
      \begin{tikzpicture}
        \begin{axis}[%
          axis lines = none,
          scatter/classes = {%
            -1={cardinal},
             1={yale}
          },
          height = 6cm,
        ]
          \addplot[scatter, only marks, mark size = 0.5pt, scatter src = explicit symbolic]
            file {Data/REDDIT-BINARY/PCA/Kernel_PCA_cleaned.txt};
        \end{axis}
      \end{tikzpicture}
    \fi
  }
  \caption{%
    Embeddings based on k-PCA for the $k_1$ kernel. \protect\subref{sfig:reddit PCA original} Every node represents a certain graph. The color
    indicates a graph from a discussion-based forum~(red) or a Q/A forum~(blue). \protect\subref{sfig:reddit PCA cleaned} We removed some
    outliers to obtain a cleaner output. It is readily visible that the two classes suffer from
    overlaps, which influence classification performance negatively.
  }
  \label{fig:reddit PCA}
\end{figure}

Here, we want to find out whether it is possible to classify the graphs using nothing but
topological information. We use the \emph{degree}, i.e., the number of neighbors, of a node in the
graph to obtain a filtration, assigning every edge the maximum of the degrees of its endpoints. We
then calculate one-dimensional persistent homology and our kernel for $p=1$ and $p=2$.
Figure~\ref{fig:reddit PCA} shows the results of applying \emph{kernel principal component
analysis}~(k-PCA)~\cite{Schoelkopf98}, which each point in the embedding corresponding to a single
graph.
A small set of outliers appears to ``skew'' the embedding~(Figure~\ref{sfig:reddit PCA original}),
but an inspection of the data shows that these graphs are extremely small~(and sparse) in contrast
to the remaining graphs. After removing them, the separation between both classes is visibly
better~(Figure~\ref{sfig:reddit PCA cleaned}).
Progress from ``left'' to ``right'' in the embedding, graphs tend to become more dense.

\begin{figure}[tbp]
  \centering
  \subfloat[$k_1$]{%
    \iffinal
      \includegraphics{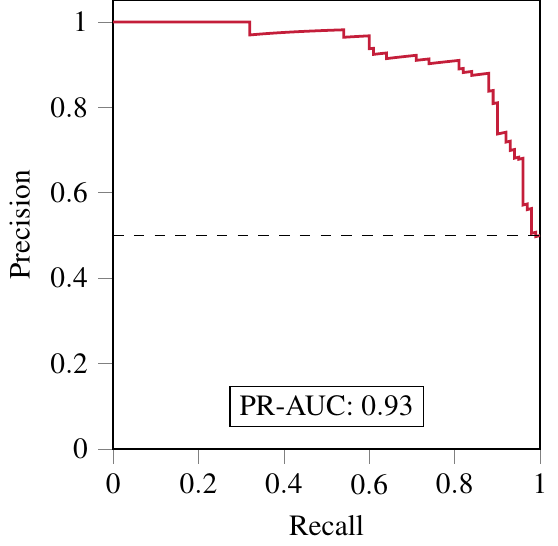}
    \else
      \begin{tikzpicture}
        \begin{axis}[%
          xmin = 0, xmax = 1.0,
          ymin = 0, ymax = 1.05,
          tick align    = outside,
          tick pos      = left,
          try min ticks = 5,
          enlargelimits     = false,
          unit vector ratio*= 1 1 1,
          xlabel = {Recall},
          ylabel = {Precision},
          scale = 0.80,
        ]
          \addplot gnuplot[thick, raw gnuplot] {%
            plot "Data/REDDIT-BINARY/PR/PR_1.txt" w l;
          };

          \draw[dashed] (0,0.5) -- (1,0.5);

          \node[fill=white, draw=black] at (0.50,0.1) {PR-AUC: $0.93$};
        \end{axis}
      \end{tikzpicture}
    \fi
  }
  \hfill
  \subfloat[$k_2$]{%
    \iffinal
      \includegraphics{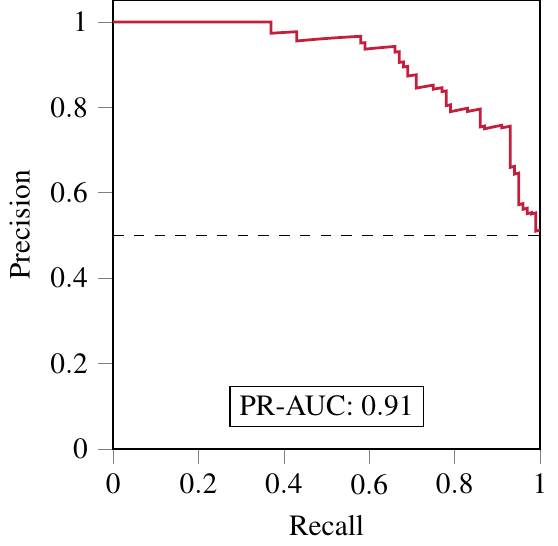}
    \else
      \begin{tikzpicture}
        \begin{axis}[%
          xmin = 0, xmax = 1.0,
          ymin = 0, ymax = 1.05,
          tick align    = outside,
          tick pos      = left,
          try min ticks = 5,
          enlargelimits     = false,
          unit vector ratio*= 1 1 1,
          xlabel = {Recall},
          ylabel = {Precision},
          scale = 0.80,
        ]
          \addplot gnuplot[thick, raw gnuplot] {%
            plot "Data/REDDIT-BINARY/PR/PR_2.txt" w l;
          };

          \draw[dashed] (0,0.5) -- (1,0.5);

          \node[fill=white, draw=black] at (0.50,0.1) {PR-AUC: $0.91$};
        \end{axis}
      \end{tikzpicture}
    \fi
  }
  \caption{%
    Precision--recall curves for both kernels on the social networks data set. Each curve also
    includes the area-under-the-curve~(AUC) value.
  }
  \label{fig:reddit precision-recall}
\end{figure}
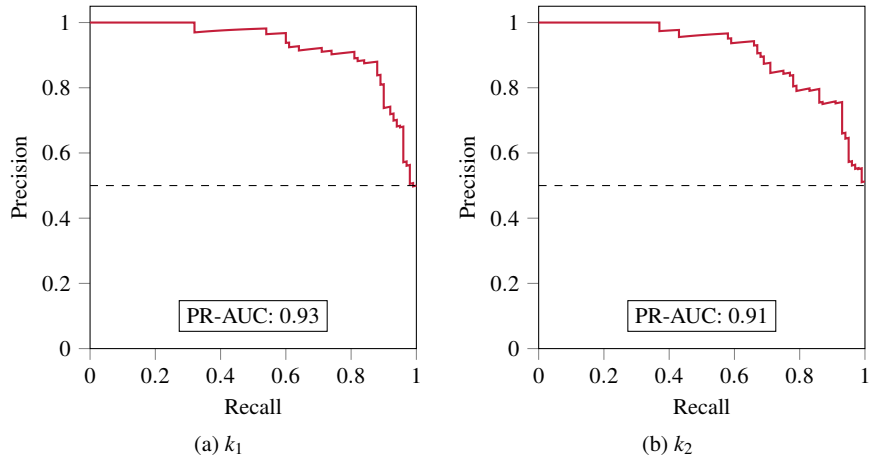

As a second application on these data, we use a kernel support vector machine to classify all
graphs, without performing outlier removal.
We split the data into training~(90\%) and test~(10\%) data, and use $4$-fold cross validation to
find the best hyperparameters. The average accuracy for $k_1$ is $0.88$, while the average accuracy
for $k_2$ is $0.81$.  \glspl{PIF} thus manage to surpass previous results by Yanardag and
Vishwanathan~\cite{Yanardag15}, which employed a computationally more expensive strategy, i.e.,
graph kernels~\cite{Sugiyama17} based on learned latent sub-structures, and obtained an average
accuracy of $0.78$ for these data.
Figure~\ref{fig:reddit precision-recall} depicts precision--recall curves for the two kernels. The
kernel $k_1$ manages to retain higher precision at higher values of recall than $k_2$, which is
again due to its lack of smoothing.

\section{Conclusion}

This paper introduced \acrfullpl{PIF}, a novel class of summarizing functions for persistence
diagrams. While being approximative by nature, we demonstrated that they exhibit beneficial
properties for data analysis, such as the possibility to perform bootstrap experiments, calculate
distances, and use kernel-based machine learning methods.
We tested the performance on various data sets and illustrated the potential of \glspl{PIF} for
topological data analysis and topological machine learning.
In the future, we want to perform a more in-depth analysis of the mathematical structure of
\glspl{PIF}, including detailed stability theorems, approximation guarantees, and a description of
their statistical properties.


\bibliographystyle{spmpsci}
\bibliography{TopoInVis2017_Persistence_Indicator_Functions}

\end{document}